\documentstyle{amsart}

\newtheorem{thm}{Theorem}
\newtheorem{prop}[thm]{Proposition}
\newtheorem{lem}[thm]{Lemma}

\theoremstyle{definition}

\numberwithin{equation}{section}

\newcommand{\bbA}{{\Bbb A}}

\newcommand{\lf}{\mathopen}

\let\r=\mathclose

\newcommand{\rank}{\operatorname{rank}}

\newcommand{\lra}{\longrightarrow}

\let\tilde=\widetilde

\let\to=\longrightarrow

\begin{document}
\title[Special groups]{On a property of special groups}
\author[Z. REICHSTEIN and B. YOUSSIN,  10-18-99]
{Z. Reichstein and B. Youssin} 
\address{Department of Mathematics, Oregon State University,
Corvallis, OR 97331\hfill\break
\hbox{{\rm\it\hskip\parindent Current mailing address\/}}: 
PMB 136, 333 South State St., Lake Oswego, OR 97034-3961, USA}
\thanks{Z. Reichstein was partially supported by NSF grant DMS-9801675 and
(during his stay at MSRI) by NSF grant DMS-9701755.}
\email{zinovy@@math.orst.edu}
\address{Department of Mathematics and Computer Science,
University of the Negev, Be'er Sheva', Israel\hfill\break
\hbox{{\rm\it\hskip\parindent Current mailing address\/}}: 
Hashofar 26/3, Ma'ale Adumim, Israel}
\email{youssin@@math.bgu.ac.il}

\begin{abstract}
Let $G$ be an algebraic group defined over an algebraically 
closed field $k$ of characteristic zero. We give a simple proof of the
following result: if $H^1(K_0, G) = \{ 1\}$ for some
finitely generated field extension $K_0/k$ of transcendence
degree $\geq 3$ then $H^1(K, G) = \{ 1 \}$ 
for every field extension $K/k$.
\end{abstract}

\maketitle

\section{Introduction}

Let $G$ be an algebraic group. J.-P. Serre stated 
the following conjectures in~\cite{serre62} (see 
also~\cite[Chapter~III]{serregc}). 

\smallskip
Conjecture I: If $G$ is connected then $H^1(K, G) = \{ 1 \}$ 
for every field $K$ of cohomological dimension $\leq 1$.

\smallskip
Conjecture II: If $G$ is semisimple, connected and simply connected 
then $H^1(K, G) = \{ 1 \}$ for every field $K$ of cohomological 
dimension $\leq 2$. 

\smallskip
Conjecture I was proved by Steinberg~\cite{steinberg0}.
Conjecture II remains open, though significant progress 
has been made in recent years; see~\cite{b-p} and~\cite{gille}.


Our main result
is a partial converse of Conjectures I and II. 
Recall that an algebraic group
$G$ is called {\em special} if $H^1(K, G) = \{ 1 \}$ for every 
field $K$. Special groups were introduced by Serre~\cite{serre1} and 
classified by Grothendieck~\cite[Section~5]{grothendieck}; 
cf.~\cite[Section~2.6]{pv}. 

\begin{thm} \label{thm1} Let $G$ be an algebraic group defined over
an algebraically closed field $k$ of characteristic zero. 
Suppose $H^1(K, G) = \{ 1 \}$ for some finitely generated 
field extension $K$ of $k$ of transcendence degree $d$.

\smallskip
(a) If $d \geq 1$ then $G$ is connected.

\smallskip
(b) If $d \geq 2$ then $G$ is simply connected.

\smallskip
(c) If $d \geq 3$ then $G$ is special.
\end{thm}

Note that the cohomological dimension of $K$ equals $d$; 
see~\cite[Section~II.4]{serregc}. Thus, informally speaking, 
the theorem may be interpreted as saying that Conjectures I and II 
cannot be extended or strengthened in a meaningful way.

Our proof of Theorem~\ref{thm1} is rather simple: the idea
is to use nontoral finite 
abelian subgroups of $G$ as obstructions to the vanishing of $H^1$. 
We remark that our argument (and, in particular, the proof of 
Lemma~\ref{lem2.1}) does not rely on canonical resolution 
of singularities; cf.~\cite[Remark~4.4]{ry}.

Ph.\ Gille recently showed us an alternative proof 
of Theorem~\ref{thm1}, based on case by case analysis and properties 
of the Rost invariant.
We would like to thank him, J.-L. Colliot-Th\'el\`ene and R. Parimala
for informative discussions.

\section{Preliminaries}

Throughout this note $k$ will denote an algebraically closed 
base field of characteristic zero. All fields, varieties, morphisms, 
algebraic groups, etc., will be assumed to be defined over $k$.

Let $G$ be an algebraic group. An abelian subgroup $A$ of $G$
is called {\em toral\/} if $A$ is contained in a torus
of $G$ and {\em nontoral\/} otherwise.

\begin{lem} \label{lem2.3} 
Let $G$ be an algebraic group, $L$ be a Levi subgroup of $G$
and $A$ be a finite abelian subgroup of $L$. If $A$ is nontoral 
in $L$ then $A$ is nontoral in $G$.
\end{lem}

\begin{pf} Assume the contrary: $A \subset T$ for some torus $T$ of
$G$. Since $T$ is reductive, it lies in a Levi subgroup $L_1$ of $G$;
see~\cite[Theorem~6.5]{ov}. 
Denote the unipotent radical of $G$ by $U$;
then $L$ and $L_1$ project isomorphically onto $G/U$.
Since $A$ is toral in $L_1$, it is toral in $G/U$, and hence, in $L$,
as claimed.
\end{pf}

Recall that a $G$-variety $X$ is an algebraic variety with 
a $G$-action; $X$ is {\em generically free} if
$G$ acts freely on a dense open subset of $X$ and 
{\em primitive} if $k(X)^G$ is a field (note that 
$X$ is allowed to be reducible). Elements of $H^1(K, G)$ are in
1---1 correspondence with $G$-{\em torsors} over $K$, i.e., 
birational classes of primitive generically free $G$-varieties 
$X$ such that $k(X)^G = K$; see e.g., \cite[Section 1.3]{popov}. 
If $X$ is a primitive generically free $G$-variety, we shall 
write $cl(X)$ for the class in $H^1(k(X)^G, G)$ given by $X$.

Our proof of Theorem~\ref{thm1} is based on the following result.

\begin{lem} (\cite[Lemma~4.3]{ry}) \label{lem2.1} 
Let $G$ be an algebraic group, $A$ be
a nontoral finite abelian subgroup of $G$ and $X$ be a generically 
free primitive $G$-variety. Suppose $A$ fixes a smooth point of $X$.
Then $cl(X) \neq 1$ in $H^1(k(X)^G, G)$.
\end{lem}

\section{Construction of a nontrivial torsor}

\begin{lem} \label{lem2.15}
Let $A$ be an abelian group of rank $r$ and let $K$ be a finitely 
generated field extension of $k$ of transcendence degree $d \geq r$. 
Then there exists an $A$-variety $Y$ such that (i) $k(Y)^A = K$ 
and (ii) $Y$ has a smooth $A$-fixed point.
\end{lem}

\begin{pf} 
Since $k$ is algebraically closed, $A$ has a faithful 
$r$-dimensional representation $V_1$.  Let $V_2$ be the trivial 
$(d-r)$-dimensional representation of $A$, and $V=V_1\oplus V_2$.
Then the (geometric)
quotient $V/A$ is isomorphic to the affine space $\bbA^d$.
Denote the origin of $V$ by $0$, and its image in $V/A$ by $\overline0$.

Let $Y_0$ be an affine variety over $k$ such that
$k(Y_0) = K$; $\dim(Y_0)=d$.
Let $y_0\in Y_0$ be a smooth point.
Identifying $V/A=\bbA^d$, we can find a dominant projection
$f \colon Y_0 \lra V/A$ such that $f(y_0)=\overline0$ and $f$ is
\'etale at $y_0$.

Now set $Y = Y_0 \times_{V/A} V$; the $A$-action 
on $Y$ is induced from $V$.  The natural projection
$Y\to Y_0$ is a rational quotient map for this action; 
see, e.g.,~\cite[Lemma~2.16(a)]{tg}. Thus $Y$ satisfies (i).
To prove (ii), set $y = (y_0, 0)$; $y$ is fixed by $A$.
The morphism $Y\to V$ is obtained from $f$ by a base
change, and hence, is \'etale at $y$; the smoothness of $V$ implies
then that $Y$ is smooth at $y$.  Thus, $y\in Y$ is a smooth point
fixed by $A$.
\end{pf}

\begin{prop} \label{prop2.2} Let $G$ be an algebraic group, $A$
be a nontoral abelian subgroup of $G$ of rank $r$, and $K/k$ be a field
extension of transcendence degree $d$. If $d \geq r$ then
$H^1(K, G) \neq \{ 1 \}$. 
\end{prop}

\begin{pf} Choose an $A$-variety $Y$ and a smooth $A$-fixed point
$y \in Y$, as in Lemma~\ref{lem2.15}.
We claim that the image of $cl(Y)$ under the natural map 
$H^1(K, A) \lra H^1(K, G)$ is nontrivial.
Indeed, recall that the image of $cl(Y)$ in $H^1(K, G)$ is $cl(X)$, 
where \[ X = G *_A Y = (G \times Y)/A \] is the (geometric) quotient 
for the $A$-action on $G \times Y$ given by $a(g, y') = (ga^{-1}, ay')$;
see~\cite[Section 4.8]{pv}.
By~\cite[Proposition~4.22]{pv}, $G*_A Y$ is smooth at $x =(1_G, y)$
since $Y$ is smooth at $y$.
Moreover, $x$ is an $A$-fixed point of $X$; thus
Lemma~\ref{lem2.1} tells us that $cl(X) \neq 1$ 
in $H^1(K, G)$, as claimed.
\end{pf}

\section{Proof of Theorem 1}

In view of Proposition~\ref{prop2.2} it is sufficient to show
that $G$ contains a nontoral finite abelian subgroup $A$, where

\smallskip
(a$'$) $\rank(A) = 1$, if $G$ is not connected,

\smallskip
(b$'$) $\rank(A) \leq 2$, if $G$ is not simply connected and

\smallskip
(c$'$) $\rank(A) \leq 3$, if $G$ is not special.

\smallskip
\noindent
Moreover, in view of Lemma~\ref{lem2.3} we only need to prove
(a$'$), (b$'$), and (c$'$) under the assumption that $G$ is reductive
(otherwise we may replace $G$ by its Levi subgroup).

\smallskip
Proof of (a$'$): Write $G = F G_0$, where
$G_0$ is the identity component of $G$ and $F$ is a finite group;
see~\cite[Proposition 7]{vinberg}. Since $G$ is disconnected, 
$F$ is not contained in $G_0$. Choose $a \in F - G_0$ and 
set $A = \lf< a\r>$.  Then $A$ is cyclic,
finite (because $a \in F$) and nontoral
(because every torus of $G$ in contained in $G_0$), as desired. 

\smallskip
Proof of (b$'$): In view of (a$'$) we may assume without loss of generality
that $G$ is connected. 
Now the desired conclusion follows from~\cite[Theorem 2.27]{steinberg}.

\smallskip
Proof of (c$'$): Suppose $G$ is not special. By 
\cite[1.5.1]{serresubgr}, $G$ has a torsion prime $p$, and
by~\cite[Theorem~2.28]{steinberg} 
$G$ has a nontoral elementary $p$-abelian subgroup $A$ of rank $\leq 3$.
see also~\cite[1.3]{serresubgr}. 
\qed

\end{document}